\documentclass{amsproc}

\newtheorem{theorem}{Theorem}[section]

\newtheorem{proposition}[theorem]{Proposition}

\theoremstyle{definition}
\newtheorem{definition}[theorem]{Definition}
\newtheorem{example}[theorem]{Example}

\theoremstyle{remark}
\newtheorem{remark}[theorem]{Remark}

\numberwithin{equation}{section}

 \DeclareMathOperator{\comp}{\#}
 \DeclareMathOperator{\id}{id}
 
 \DeclareMathOperator{\tri}{\bigtriangledown}
\newcommand{\F}{\mathbf{F}}
\newcommand{\Or}{\mathcal{O}}
\newcommand{\Z}{\mathbf{Z}}

\begin{document}

\title{Orientals}

\author{Richard Steiner}
\address{Department of Mathematics\\University of Glasgow\\University Gardens\\
Glasgow\\Scotland G12 8QW} \email{r.steiner@maths.gla.ac.uk}

\subjclass[2000]{Primary 18D05}

\dedicatory{For Ross Street on his sixtieth birthday}

\date{}

\begin{abstract} The orientals or oriented simplexes are a family
of strict omega-categories constructed by Ross Street. We show
that the category of orientals is isomorphic to a subcategory of
the category of chain complexes. This leads to a very simple
combinatorial description of the morphisms between orientals. We
also show that the category of orientals is the closure of the
category of simplexes under certain filler operations which
represent complicial operations.
\end{abstract}

\maketitle

\section{Introduction}
\label{Introduction}

The orientals or oriented simplexes are a family of strict
$\omega$-categories $\Or_0$, $\Or_1$, \dots\ constructed by
Street~\cite{StrSimp}. In \cite{Ste} it is shown that there is a
full subcategory of the category of strict $\omega$-categories
which is isomorphic to a category of chain complexes with
additional structure, and it is shown that the orientals lie in
this full subcategory. Implicitly, this describes the morphisms
between orientals. The object of the present paper is to make the
description of the morphisms more explicit. It turns out that they
have a very simple combinatorial description; see
Theorem~\ref{oriental morphisms} below. This supports Street's
assertion in~\cite{StrSimp} that the orientals are fundamental
structures of nature.

Let $\Or$ be the category of orientals, let $\Delta$ be the
simplex category, and let $\Z\Delta$ be the category whose
morphisms are integer linear combinations of morphisms
in~$\Delta$; the morphisms in $\Z\Delta$ are to be composed
bilinearly. Theorem~\ref{oriental morphisms} shows that $\Or$~is a
subcategory of $\Z\Delta$ containing~$\Delta$. The morphisms
in~$\Or$ can therefore be expressed in terms of those in~$\Delta$
by using the operations of addition and subtraction, but these
operations are defined only in the larger category $\Z\Delta$. We
will however show that $\Or$~is the closure of~$\Delta$ under a
family of operations internal to~$\Or$. These operations are
called filler operations and they are the universal examples for
Verity's complicial operations~\cite{V}. Verity shows that
$\omega$-category nerve structures on simplicial sets are
equivalent to families of complicial operations with suitable
properties, and our closure result was inspired by Verity's
theorem.

In Section~\ref{General} we recall the relevant results relating
chain complexes and $\omega$-categories; in Section~\ref{objects}
we show how the individual orientals can be expressed in terms of
chain complexes; and in Section~\ref{morphisms} we describe the
morphisms between the orientals. Finally, Section~\ref{fillers and
morphisms} contains the material on fillers.

\section{Chain complexes and omega-categories}
\label{General}

In~\cite{Ste} there is a description of a functor~$\nu$ from a
category of chain complexes with additional structure to the
category of strict $\omega$-categories; it is modelled on Street's
theory of parity complexes~\cite{StrPar}. In this section we
recall some results from~\cite{Ste} concerning the functor~$\nu$.
All $\omega$-categories in this paper are strict
$\omega$-categories.

We restrict attention to augmented chain complexes concentrated in
nonnegative dimensions whose chain groups are free abelian groups
with prescribed bases. Throughout this section, let
$K=(K,\partial,\epsilon)$ be such a chain complex. We will say
that $K$~itself has a basis, consisting of the disjoint union of
the bases for the individual chain groups. We make each chain
group~$K_q$ into a partially ordered abelian group as follows: if
$x\in K_q$, then $x\geq 0$ if and only if $x$~is a sum of basis
elements. We will now define the associated $\omega$-category $\nu
K$, using Street's one-sorted notion of $\omega$-category
\cite[1]{StrSimp}; thus $\nu K$ is a single set which serves as
the set of morphisms for every category in a family indexed by the
nonnegative integers.

\begin{definition}
\label{nu} The members of the $\omega$-category $\nu K$ are the
double sequences
$$x=(x^-_0,x^+_0\mid x^-_1,x^+_1\mid \ldots\,)$$
satisfying the following conditions:
\begin{enumerate}
 \item $x^-_q,x^+_q\in K_q$;
 \item there exists~$N$ such that $x^-_q=x^+_q=0$ for $q>N$;
 \item $x^+_q-x^-_q=\partial x^-_{q+1}=\partial x^+_{q+1}$ for
 $q\geq 0$;
 \item $x^-_q\geq 0$ and $x^+_q\geq 0$ for $q\geq 0$;
 \item $\epsilon x^-_0=\epsilon x^+_0=1.$
\end{enumerate}
For $p\geq 0$ we make $\nu K$ into the morphism set of a category
as follows: the left identity $d^-_p x$ and the right identity
$d^+_p x$ of the element $x=(x^-_0,x^+_0\mid
x^-_1,x^+_1\mid\ldots\,)$ are given by
$$d^-_p x
 =(x^-_0,x^+_0\mid\ldots
 \mid x^-_{p-1},x^+_{p-1}\mid x^-_p,x^-_p\mid 0,0\mid\ldots\,)$$
and
$$d^+_p x
 =(x^-_0,x^+_0\mid\ldots
 \mid x^-_{p-1},x^+_{p-1}\mid x^+_p,x^+_p\mid 0,0\mid\ldots\,);$$
if $x$~and~$y$ are elements such that $d^+_p x=d^-_p y$, say
$d^+_p x=d^-_p y=w$, then the composite $x\comp_p y$ is given by
$$x\comp_p y=x-w+y,$$
where the addition and subtraction are performed termwise.
\end{definition}

In the multi-sorted description of $\nu K$ the $p$-cells are the
elements of the form
$$(x^-_0,x^+_0,\mid\ldots\mid x^-_{p-1},x^+_{p-1}\mid x_p,x_p\mid
 0,0\mid\ldots\,).$$

\begin{remark} Definition~\ref{nu} is based on the well-known
equivalence between chain complexes concentrated in nonnegative
dimensions and $\omega$-category objects in the category of
abelian groups; see \cite{BH} for example. Indeed, if $K$~is a
chain complex concentrated in nonnegative degrees, then the
equivalent $\omega$-category object is defined in the same way as
$\nu K$, except that conditions (4) and~(5) are omitted.
\end{remark}

In order to construct some members of $\nu K$, we now introduce
some notation. Given a chain~$x$, we write $\partial^+ x$ and
$\partial^- x$ for the positive and negative parts of the boundary
$\partial x$; thus $\partial^+ x$ and $\partial^- x$ are sums of
basis elements with no common terms such that
$$\partial x=\partial^+ x-\partial^- x.$$
Given a $p$-dimensional chain~$x$, we write $\langle x\rangle$ for
the double sequence given by
$$\langle x\rangle
 =\bigl((\partial^-)^p x,(\partial^+)^p x\mid
 \ldots\mid
 \partial^- x,\partial^+ x\mid
 x,x\mid
 0,0\mid
 \ldots\,\bigr),$$
and we observe that $\langle x\rangle$ is a member of $\nu K$ if
and only if $\epsilon(\partial^-)^p x=\epsilon(\partial^+)^p x=1$.
This observation motivates the following definition; recall that
the prescribed basis for~$K$ means the disjoint union of the
prescribed basis for the individual chain groups.

\begin{definition} \label{unital} The prescribed basis for the
chain complex~$K$ is \emph{unital} if $\epsilon(\partial^-)^p
b=\epsilon(\partial^+)^p b=1$ whenever $b$~is a basis element of
dimension~$p$. If $K$~has a unital basis, then the
elements~$\langle b\rangle$ of $\nu K$ corresponding to the basis
elements~$b$ are called \emph{atoms}.
\end{definition}

The main result says that $\nu$~restricts to an equivalence if one
imposes suitable restrictions on the bases of the chain complexes
involved. It suffices that the bases should be both unital and
strongly loop-free, where strong loop-freeness is defined as
follows.

\begin{definition} \label{loop-free} The prescribed basis for the
chain complex~$K$ is \emph{strongly loop-free} if it has a partial
ordering such that $a<b$ whenever $a$~is an element with a
negative coefficient in $\partial b$ or $b$~is an element with a
positive coefficient in $\partial a$.
\end{definition}

Note that the partial ordering in Definition~\ref{loop-free} is
quite different from the partial orderings of the individual chain
groups.

By combining these restrictions we obtain a category~$\F$ of chain
complexes as follows.

\begin{definition} \label{category F} The objects of~$\F$ are the
augmented chain complexes of abelian groups concentrated in
nonnegative degrees together with prescribed strong\-ly loop-free
unital bases. The morphisms of~$\F$ are the augmentation and
order-preserving chain maps.
\end{definition}

The main result \cite[5.11]{Ste} is now as follows.

\begin{theorem} \label{main} The construction~$\nu$ is a full and
faithful embedding of the category~$\F$ in the category of strict
$\omega$-categories.
\end{theorem}

We will also use the following result \cite[5.6]{Ste}.

\begin{theorem} \label{atom generation} If $K$~is an object of~$\F$,
then the $\omega$-category $\nu K$ is generated under composition
by its atoms.
\end{theorem}

\section{The construction of the orientals} \label{objects}

Let $K[n]$ be the chain complex of the standard $n$-simplex, where
$n$~is a nonnegative integer. We will show that $K[n]$ is an
object in the category~$\F$ and that $\nu K[n]$ is isomorphic to
the $n$th oriental~$\Or_n$. Essentially, we are working with free
abelian groups whose bases are certain sets constructed by Street
in~\cite{StrSimp}. The abelian group context makes the theory
rather easier, but the proof that $K[n]$ is in~$\F$ uses the same
computations as in~\cite{StrSimp}.

We recall that $K[n]$ is an augmented chain complex of free
abelian groups; the prescribed basis of $K[n]_q$ consists of the
ordered $(q+1)$-tuples of integers $[a_0,\ldots,a_q]$ such that
$0\leq a_0<a_1<\ldots<a_q\leq n$; the boundary $\partial\colon
K[n]_q\to K[n]_{q-1}$ for $q>0$ is given by the alternating sum
formula
$$\partial[a_0,\ldots,a_q]
 =[a_1,\ldots,a_q]
 -[a_0,a_2,\ldots,a_q]
 +\ldots
 +(-1)^q[a_0,\ldots,a_{q-1}];$$
the augmentation is given by
$$\epsilon[a_0]=1.$$

\begin{theorem} \label{simplex bases} The prescribed basis for
$K[n]$ is strongly loop-free and unital, so that $K[n]$ is an
object of~$\F$.
\end{theorem}

\begin{proof} To show that the basis is strongly loop-free, we
must find a partial ordering of the basis elements such that $a<b$
if $a$~is a term in $\partial b$ with a negative coefficient or if
$b$~is a term in $\partial a$ with a positive coefficient. There
is in fact a total ordering with these properties; it can be
described recursively as follows. Let $a=[a_0,\ldots,a_p]$ and
$b=[b_0,\ldots,b_q]$ be distinct basis elements. Then $a<b$ if
$a_0<b_0$, or if $a_0=b_0$ and $p=0$, or if $a_0=b_0$ and $p,q>0$
and $[a_1,\ldots,a_p]>[b_1,\ldots,b_q]$.

To show that the basis is unital, let $a=[a_0,\ldots,a_p]$ be a
basis element; we must show that $\epsilon(\partial^-)^p
a=\epsilon(\partial^+)^p a=1$. In order to do this, we compute
$(\partial^-)^q a$ and $(\partial^+)^q a$ for $0<q\leq p$
inductively. We find that $(\partial^-)^q a$ and $(\partial)^q a$
are the sums of terms $[a_{i(0)},\ldots,a_{i(p-q)}]$ such that the
indices of the omitted elements $a_{j(1)},\ldots,a_{j(q)}$ form an
increasing sequence $j(1)<\ldots<j(q)$ of integers of alternating
parity; to be precise, we get the terms of $(\partial^-)^q a$ by
taking $j(1)$ to be odd, $j(2)$ to be even, etc., and we get the
terms of $(\partial^+)^q a$ by taking $j(1)$ to be even, $j(2)$ to
be odd, etc. In particular, if $p>0$ then we get $(\partial^-)^p
a=[a_0]$ and $(\partial^+)^p a=[a_p]$, and these formulae are
obviously valid for $p=0$ as well. It follows that
$\epsilon(\partial^-)^p a=\epsilon(\partial^+)^p a=1$ as required.
\end{proof}

In \cite[6]{Ste} an axiomatic characterisation of $\nu K[n]$ is
used to show that it is  isomorphic to the $n$th oriental. Here we
sketch a direct comparison.

\begin{theorem} \label{orientals} The $\omega$-category $\nu K[n]$ is
isomorphic to the $n$th oriental~$\Or_n$.
\end{theorem}

\begin{proof} The elements of~$\Or_n$ (see \cite[2]{StrSimp})
are double sequences
$$a=(a^1_0,a^0_0\mid a^1_1,a^0_1\mid\ldots\,),$$
where $a^1_q$~and~$a^0_q$ are certain finite sets of
$q$-dimensional basis elements for $K[n]$. Associated to such a
double sequence~$a$ there is a double sequence
$$x=(x^-_0,x^+_0\mid x^-_1,x^+_1\mid\ldots\,)$$
such that $x^-_q$~and~$x^+_q$ are $q$-dimensional chains in
$K[n]$; one takes $x^-_q$~and~$x^+_q$  to be the sums of the
members of~$a^1_q$ and of~$a^0_q$ respectively. From the
definition of~$\Or_n$, one can check that $x$~is in $\nu K[n]$, so
we have obtained a function $i\colon\Or_n\to\nu K[n]$, and it is
clear that $i$~is injective. By comparing the definitions, one
sees that $i$~is a morphism of $\omega$-categories. One can also
check that the atoms are in the image of~$i$. Since $\nu K[n]$ is
generated by its atoms (see Theorem~\ref{atom generation}), it
follows that $i$~is an isomorphism. Therefore $\nu K[n]\cong\Or_n$
as required.
\end{proof}

\section{Morphisms between orientals} \label{morphisms}

By combining Theorems \ref{orientals} and~\ref{main}, we see that
the category of orientals is isomorphic to the category of
augmentation and order-preserving chain maps between the chain
complexes $K[n]$. In this section we get a simple description of
the chain maps, and from this we get a simple description of the
morphisms between the orientals.

We begin with the simplex category~$\Delta$. The objects
of~$\Delta$ are the nonnegative integers $0,1,\ldots\,$, the
morphism set $\Delta(m,n)$ is the set of non-decreasing functions
$$f\colon\{0,1,\ldots,m\}\to\{0,1,\ldots,n\},$$
and the morphisms are composed in the obvious way. For the
morphism in $\Delta(m,n)$ given by $i\mapsto f_i$ we will use the
notation $(f_0,\ldots,f_m)$; thus the morphisms in $\Delta(m,n)$
are the ordered $(m+1)$-tuples of integers $(f_0,\ldots,f_m)$ such
that
$$0\leq f_0\leq f_1\leq\ldots\leq f_m\leq n.$$

There is a standard functor from~$\Delta$ to chain complexes
sending~$n$ to $K[n]$, as follows. Let $f=(f_0,\ldots,f_m)$ be a
morphism in $\Delta(m,n)$ and let $[a(0),\ldots,a(q)]$ be a basis
element for $K[m]$. If $f_{a(0)}<f_{a(1)}<\ldots<f_{a(q)}$ then
$$f[a(0),\ldots,a(q)]=[f_{a(0)},\ldots,f_{a(q)}];$$
if $f_{a(0)},\ldots,f_{a(q)}$ are not distinct then
$$f[a(0),\ldots,a(q)]=0.$$

Now let $\Z\Delta$ be the category with objects $0,1,2,\ldots\,$
such that the set of morphisms $\Z\Delta(m,n)$ is the free abelian
group with basis $\Delta(m,n)$. The morphisms in $\Z\Delta$ are
composed bilinearly; that is,
$$\Bigl(\sum_g\mu_g g\Bigr)\circ\Bigl(\sum_f\lambda_f f\Bigr)
 =\sum_h\Bigl(\sum_{g\circ f=h}\mu_g\lambda_f\Bigr)h.$$
Since the chain maps from $K[m]$ to $K[n]$ form an abelian group,
the functor $n\mapsto K[n]$ from~$\Delta$ to chain complexes
extends to a functor on $\Z\Delta$; in other words, there is a
natural homomorphism from $\Z\Delta(m,n)$ to the group of chain
maps from $K[m]$ to $K[n]$. We will now show that this
homomorphism is an isomorphism.

\begin{theorem} \label{chain maps} The group of chain maps from
$K[m]$ to $K[n]$ is naturally isomorphic to $\Z\Delta(m,n)$.
\end{theorem}

\begin{proof} Let $A_q$~and~$B_q$ be the bases for $K[m]_q$ and
$K[n]_q$ respectively, let $A'_q$ be the subset of~$A_q$
consisting of the elements $[a_0,\ldots,a_q]$ such that $a_0=0$,
and let $K[m]'_q$ be the subgroup of $K[m]_q$ generated by the
members of $A'_q$. If $[a_0,\ldots,a_q]$  is a basis element of
$K[m]_q$ such that $a_0>0$, then
$$[a_0,\ldots,a_q]=\partial[0,a_0,\ldots,a_q]+s$$
with $s\in K[m]'_q$, so a chain map from $K[m]$ to $K[n]$ is
uniquely determined by its restrictions to the subgroups
$K[m]'_q$. It is therefore sufficient to show that these
restrictions induce an isomorphism
$$\Z\Delta(m,n)\to\bigoplus_q\hom(K[m]'_q,K[n]_q).$$

To do this, let $f$ be a morphism in $\Delta(m,n)$ and let the
size of the image of~$f$ be $q+1$. To~$f$ we associate basis
elements $a=[a_0,\ldots,a_q]$ in~$A'_q$ and $b=[b_0,\ldots,b_q]$
in~$B_q$ as follows: $b_0,\ldots,b_q$ are the members of the image
of~$f$ in ascending order, and $a_i$~is the smallest member of the
inverse image $f^{-1}(b_i)$. We find that these assignments
produce a bijection
$$\Delta(m,n)\to\coprod_q(A'_q\times B_q);$$
we also find that $f[a_0,\ldots,a_q]=[b_0,\ldots,b_q]$, that
$f[i_0,\ldots,i_r]=0$ if $[i_0,\ldots,i_r]$ is a member of~$A'_r$
for some $r>q$, and that $f[i_0,\ldots,i_q]=0$ if
$[i_0,\ldots,i_q]$ is a member of~$A'$ which precedes
$[a_0,\ldots,a_q]$ lexicographically (i.e.\ there exists~$r$ such
that $i_0=a_0$, \dots, $i_{r-1}=a_{r-1}$, $i_r<a_r$). It follows
that the morphism
$$\Z\Delta(m,n)\to\bigoplus_q\hom(K[m]'_q,K[n]_q)$$
is an isomorphism, as required. This completes the proof.
\end{proof}

Theorem~\ref{chain maps} shows that there is a category $\Z\Delta$
with objects $0,1,2,\ldots\,$ and morphism sets $\Z\Delta(m,n)$.
According to Theorem~\ref{main}, the category of orientals is
isomorphic to the subcategory of $\Z\Delta$ consisting of the
augmentation and order-preserving morphisms. It is easy to
characterise these morphisms, and the result is as follows.

\begin{theorem} \label{oriental morphisms} The set of morphisms
from the oriental~$\Or_m$ to the oriental~$\Or_n$ is naturally
isomorphic to the subset of $\Z\Delta(m,n)$ consisting of the
linear combinations~$x$ such that
\begin{enumerate}
 \item the sum of the coefficients in~$x$ is~$1$\textup{;}
 \item for all injective morphisms~$f$ in~$\Delta$ with codomain~$m$
 the coefficients of the injective morphisms in the composite $x\circ
 f$ are nonnegative.
\end{enumerate}
\end{theorem}

The subset of $\Z\Delta(m,n)$ described in Theorem~\ref{oriental
morphisms} will be denoted $\Or(m,n)$, and the category of
orientals will be denoted~$\Or$; thus $\Or$~has objects~$\Or_n$
and morphism sets $\Or(m,n)$.

\begin{example} \label{filler examples}
If $q>0$ and $0\leq i_0<i_1<\ldots<i_q\leq n$ then there is a
morphism~$x$ in $\Or(1,n)$ given by
$$x=(i_0,i_1)-(i_1,i_1)+(i_1,i_2)-(i_2,i_2)+\ldots+(i_{q-1},i_q);$$
indeed the sum of the coefficients is~$1$, there are no injective
morphisms with negative coefficients in $x\circ(0,1)=x$, and there
are no injective morphisms with negative coefficients in
$x\circ(0)=(i_0)$ or $x\circ(1)=(i_q)$ either.

Similarly there is a morphism~$y$ in $\Or(2,2)$ given by
$$y=(0,1,1)-(1,1,1)+(1,1,2).$$
\end{example}

\section{Fillers}
\label{fillers and morphisms}

According to Theorem~\ref{oriental morphisms}, the morphisms in
the category~$\Or$ of orientals are integer linear combinations of
morphisms in the subcategory~$\Delta$; but one cannot define
addition and subtraction inside~$\Or$, only in the larger category
$\Z\Delta$. In this section we will show how to generate~$\Or$
from~$\Delta$ by using operations defined inside~$\Or$.

We first observe that the sets $\Z\Delta(m,n)$ for a fixed~$n$ and
varying~$m$ form a simplicial set, and that the $\Or(m,n)$ form a
simplicial subset. That is to say, there are face operations
$\partial_i\colon\Z\Delta(m,n)\to\Z\Delta(m-1,n)$ for $m>0$ and
$0\leq i\leq m$ given by
$$\partial_i x=x\circ(0,\ldots,i-1,i+1,\ldots,m),$$
there are degeneracy operations
$\epsilon_i\colon\Z\Delta(m,n)\to\Z\Delta(m+1,n)$ for $0\leq i\leq
m$ given by
$$\epsilon_i x=x\circ(0,\ldots,i-1,i,i,i+1,\ldots,m),$$
the subcategory~$\Or$ is closed under these operations, and the
usual simplicial identities apply. In this paper we will use the
following simplicial identities:
$$\partial_i\epsilon_{i+1}=\epsilon_i\partial_i,\qquad
 \partial_i\epsilon_i=\partial_{i+1}\epsilon_i=\id,\qquad
 \partial_{i+2}\epsilon_i=\epsilon_i\partial_{i+1}.$$
We now define two further families of operations.

\begin{definition} \label{filler composite definition} Let
$x$~and~$y$ be morphisms in $\Z\Delta(m,n)$ such that $\partial_i
x=\partial_{i+1}y$ for some~$i$ with $0\leq i\leq m-1$. Then
$x\tri_i y$ is the morphism in $\Z\Delta(m+1,n)$ given by
$$x\tri_i y
 =\epsilon_{i+1}x-\epsilon_i\epsilon_i\partial_i x+\epsilon_i y$$
and $x\vee_i y$ is the morphism in $\Z\Delta(m,n)$ given by
$$x\vee_i y=x-\epsilon_i\partial_i x+y.$$
We call $x\tri_i y$ a \emph{filler} and $x\vee_i y$ a
\emph{pasting}.
\end{definition}

For example, let $x$~and~$y$ be the morphisms of
Example~\ref{filler examples}; then
$$x=(i_0,i_1)\vee_0\ldots\vee_0(i_{q-1},i_q)$$
and $y=(0,1)\tri_0(1,2)$. (No brackets are needed in the
expression for~$x$, since one finds that $\vee_0$~is associative.)
As models for use in later arguments, we also give the following
examples.

\begin{example} \label{useful}
If $a\leq b$ then
\begin{align*}
 &(a,a)\tri_0(a,b)=(a,a,b),\\
 &(a,b)\tri_0(b,b)=(a,b,b),\\
 &(a,a)\vee_0(a,b)=(a,b),\\
 &(a,b)\vee_0(b,b)=(a,b).
\end{align*}
\end{example}

From the simplicial identities we deduce the following result.

\begin{proposition} \label{faces} If $x$~and~$y$ are morphisms in
$\Z\Delta(m,n)$ such that $\partial_i x=\partial_{i+1}y$, then
\begin{align*}
 &\partial_i(x\tri_i y)=y,\\
 &\partial_{i+1}(x\tri_i y)=x\vee_i y,\\
 &\partial_{i+2}(x\tri_i y)=x.
\end{align*}
\end{proposition}

The object of this section is to show that the category~$\Or$ is
the closure of its subcategory~$\Delta$ under the
operations~$\tri_i$. We begin by verifying that $\Or$~is closed
under the filler and pasting operations.

\begin{proposition} \label{closure} If $x$~and~$y$ are morphisms
in $\Or(m,n)$ such that $\partial_i x=\partial_{i+1}y$ then
$x\tri_i y\in\Or(m+1,n)$ and $x\vee_i y\in\Or(m,n)$.
\end{proposition}

\begin{proof} Since $x\vee_i y=\partial_{i+1}(x\tri_i y)$, it
suffices to consider $x\tri_i y$. It is clear that the sum of the
coefficients in $x\tri_i y$ is equal to~$1$, and it remains to
show that injective morphisms have nonnegative coefficients in
$(x\tri_i y)\circ f$ when $f$~is an injective morphism in~$\Delta$
with codomain $m+1$. If $i+2$ is not in the image of~$f$, then
this holds because $\partial_{i+2}(x\tri_i y)=x\in\Or(m,n)$. If
$i$~is not in the image of~$f$, then this holds similarly because
$\partial_i(x\tri_i y)=y$. If both $i$~and $i+2$ occur in the
image of~$f$ but $i+1$ does not, then the injective terms in
$(x\tri_i y)\circ f$ are obtained by adding those in
$(\epsilon_{i+1}x)\circ f$ and $(\epsilon_i y)\circ f$, so again
they have nonnegative coefficients. Finally, if $i$, $i+1$ and
$i+2$ are all in the image of~$f$, then $(x\tri_i y)\circ f$ has
no injective terms at all. This completes the proof.
\end{proof}

Conversely, we now give a sequence of results designed to express
the morphisms of~$\Or$ in terms of fillers and pastings.

\begin{proposition} \label{tails} Let $x$ be a morphism in
$\Z\Delta(m,n)$ with $m>0$ such that injective morphisms have
nonnegative coefficients in $x\circ f$ for all injective
morphisms~$f$ in~$\Delta$ with codomain~$m$. In the obvious
notation, write
$$x=(x_0,0)+\ldots+(x_n,n)$$
with $x_i\in\Z\Delta(m-1,i)$. Then injective morphisms have
nonnegative coefficients in
$$(x_r+x_{r+1}+\ldots+x_n)\circ g$$
for $0\leq r\leq n$ and for all injective morphisms~$g$
in~$\Delta$ with codomain $m-1$.
\end{proposition}

\begin{proof} Let $g=(g_0,\ldots,g_q)$, and let $a=(a_0,\ldots,a_q)$
be an injective morphism in $\Delta(q,n)$. If $a_q<r$ then the
coefficient of~$a$ in $(x_r+\ldots+x_n)\circ g$ is the sum over
$r\leq i\leq n$ of the coefficients of $(a_0,\ldots,a_q,i)$ in
$x\circ(g_0,\ldots,g_q,m)$; if $a_q\geq r$ then the coefficient
of~$a$ in $(x_r+\ldots+x_n)\circ g$ is the same as its coefficient
in $x\circ(g_0,\ldots,g_q)$. The result follows.
\end{proof}

\begin{proposition} \label{zero sums} Let $x$ be a morphism in
$\Z\Delta(m,n)$ such that injective morphisms have nonnegative
coefficients in $x\circ f$ for all injective morphisms~$f$
in~$\Delta$ with codomain~$m$ and such that the sum of the
coefficients in~$x$ is zero. Then $x=0$.
\end{proposition}

\begin{proof}
We use induction on~$m$. Note first that every morphism in
$\Delta(0,n)$ is injective, so $x\circ(m)$ is a linear combination
with nonnegative coefficients such that the sum of the
coefficients is zero. It follows that $x\circ(m)=0$.

Suppose that $m=0$. Then $x=x\circ(m)$, so $x=0$.

Now suppose that $m>0$. Write~$x$ in the form
$$x=(x_0,0)+\ldots+(x_n,n)$$
as in Proposition~\ref{tails}. Since $x\circ(m)=0$, the sum of the
coefficients in~$x_i$ is zero for $0\leq i\leq m$. From
Proposition~\ref{tails} and the inductive hypothesis, it follows
that $x_r+x_{r+1}+\ldots+x_n=0$ for $0\leq r\leq n$, and it then
follows that $x_i=0$ for $0\leq i\leq n$. Therefore $x=0$, as
required.

This completes the proof.
\end{proof}

\begin{proposition} \label{first and last} Let $x$ be a morphism
in $\Or(m,n)$, let $s$ be the smallest integer to appear in the
terms of~$x$, and let $t$ be the largest integer to appear in the
terms of~$x$. Then $x[0]=[s]$ and $x[m]=[t]$.
\end{proposition}

\begin{proof} Since every morphism in $\Delta(0,n)$ is injective,
$x\circ(0)$ and $x\circ(m)$ are linear combinations with
nonnegative integer coefficients. Since the sums of the
coefficients are~$1$, they reduce to single terms, so that
$x\circ(0)=(a)$ and $x\circ(m)=(b)$ for some integers $a$~and~$b$.
It follows that $x[0]=[a]$ and $x[m]=[b]$, and it now suffices to
show that $s=a$ and $t=b$. We will show that $t=b$; the proof that
$s=a$ is similar.

Suppose first that $m=0$. Then $x=x\circ(m)=(b)$, so $t=b$.

Now suppose that $m>0$. Write~$x$ in the form
$$x=(x_0,0)+\ldots+(x_n,n)$$
as in Proposition~\ref{tails}. Since $t$~is the largest integer to
appear in the terms of~$x$, we have $x_t\neq 0$ and $x_i=0$ for
$t<i\leq n$. It follows from Proposition~\ref{tails} that the
injective terms in $x_t\circ f$ have nonnegative coefficients for
all injective morphisms~$f$ in~$\Delta$ with codomain $m-1$. Since
$x_t\neq 0$, it follows from Proposition~\ref{zero sums} that the
sum of the coefficients of~$x_t$ is not zero; in other words $[t]$
has a non-zero coefficient in $x[m]$. Since $x[m]=[b]$, we must
have $t=b$ as required. This completes the proof.
\end{proof}

\begin{proposition} \label{constant}
Let $x$ be a morphism in $\Or(m,n)$ with $x[m]=[t]$ such that
$(t,\ldots,t)$ has a non-zero coefficient in~$x$. Then
$x=(t,\ldots,t)$.
\end{proposition}

\begin{proof} By Proposition~\ref{first and last}, $t$~is the largest
integer to appear in the terms of~$x$. The coefficient of~$[t]$ in
$x[0]$ is therefore the coefficient of $(t,\ldots,t)$ in~$x$,
which we are assuming to be non-zero. By Proposition~\ref{first
and last} again, $t$~is also the smallest integer to appear in the
terms of~$x$. It follows that $x=(t,\ldots,t)$, as required.
\end{proof}

\begin{proposition} \label{factor filler start}
Let $x$ be a morphism in $\Or(m,n)$ with $x[m]=[t]$, and let $r$
be an integer with $0\leq r\leq m-2$ such that $x$~is a linear
combination of terms $(a_0,\ldots,a_m)$ with $a_r<t$. Then there
is a factorisation
$$x=u\vee_r v=u\vee_r(\partial_{r+2}v\tri_r\partial_r v)$$
in~$\Or$ such that $u[m]=[t]$ and $u$~is a linear combination of
terms $(a_0,\ldots,a_m)$ with $a_{r+1}<t$.
\end{proposition}

\begin{proof} Let $\alpha,\beta\colon\Z\Delta(m,n)\to\Z\Delta(m,n)$ be the
linear functions given on morphisms $(a_0,\ldots,a_m)$ in
$\Delta(m,n)$ as follows:
\begin{align*}
&\alpha(a_0,\ldots,a_m)=
 \left\{\begin{alignedat}{2}
 &(a_0,\ldots,a_{r-1},a_r,a_{r+1},a_{r+2}\ldots,a_m)
 &\quad &\text{if $a_{r+1}<t$,}\\
 &(a_0,\ldots,a_{r-1},a_r,a_r,a_{r+2}\ldots,a_m)
 &&\text{if $a_{r+1}\geq t$,}
 \end{alignedat}\right.\\
&\beta(a_0,\ldots,a_m)=
 \left\{\begin{alignedat}{2}
 &(a_0,\ldots,a_{r-1},a_{r+1},a_{r+1},a_{r+2}\ldots,a_m)
 &\quad &\text{if $a_{r+1}<t$,}\\
 &(a_0,\ldots,a_{r-1},a_r,a_{r+1},a_{r+2}\ldots,a_m)
 &&\text{if $a_{r+1}\geq t$.}
 \end{alignedat}\right.
\end{align*}
For a morphism $a=(a_0,\ldots,a_m)$ in $\Delta(m,n)$ with $a_r<t$
and $a_m\leq t$ we find that
$$a
 =\alpha a\vee_r\beta a
 =\alpha a\vee_r(\partial_{r+2}\beta a\tri_r\partial_r\beta a)$$
(the calculation is as in Example~\ref{useful}). Since $x$~is a
linear combination of morphisms of this type, we get a
factorisation
$$x=u\vee_r v=u\vee_r(\partial_{r+2}v\tri_r\partial_r v)$$
in $\Z\Delta(m,n)$, where $u=\alpha x$ and $v=\beta x$. Clearly
$u[m]=x[m]$, so that $u[m]=[t]$, and it remains to show that the
factors are in~$\Or$. Since $\Or$~is closed under the
operations~$\partial_i$, it suffices to show that $u$~and~$v$ are
in $\Or(m,n)$ by verifying the conditions of Theorem~\ref{oriental
morphisms}. Now, it is clear that the sums of the coefficients in
$u$~and~$v$ are equal to~$1$, so it remains to consider $u\circ f$
and $v\circ f$ when $f=(f_0,\ldots,f_q)$ is injective with
codomain~$m$; we must show that every injective morphism
$i=(i_0,\ldots,i_q)$ with codomain~$n$ has nonnegative
coefficients in $u\circ f$ and $v\circ f$.

If $r+1$ is not in the image of~$f$, then $u\circ f=x\circ f$, so
the result holds for $u\circ f$. Similarly, if $r$~is not in the
image of~$f$, then $v\circ f=x\circ f$, so the result holds for
$v\circ f$.

Now suppose that both $r$~and $r+1$ are in the image of~$f$, say
$f_{j-1}=r$ and $f_j=r+1$. If $i_j<t$, then the coefficient of~$i$
in $u\circ f$ is as in $x\circ f$ and the coefficient of~$i$ in
$v\circ f$ is zero. If $i_j\geq t$ then the coefficient of~$i$ in
$u\circ f$ is zero and the coefficient of~$i$ in $v\circ f$ is as
in $x\circ f$. In both cases, the coefficients of~$i$ in $u\circ
f$ and $v\circ f$ are therefore nonnegative.

Next we consider $u\circ f$ when the image of~$f$ contains $r+1$
but not~$r$, say $f_j=r+1$. If $i_j\geq t$, then the coefficient
of~$i$ in $u\circ f$ is zero; if $i_j<t$ and $i_k<t$ for some
$k>j$, then the coefficient of~$i$ in $u\circ f$ is the same as in
$x\circ f$; if $j=q-1$ and $i_j<t$ and $i_q=t$, then the
coefficient of~$i$ in $u\circ f$ is the sum of the coefficients
of~$i$ in $x\circ f$ and in $x\circ(f_0,\ldots,f_{q-2},r,r+1)$; if
$j=q$ and $i_j<t$, then the coefficient of~$i$ in $u\circ f$ is
the sum of the coefficients of~$i$ in $x\circ f$ and of
$(i_0,\ldots,i_q,t)$ in $x\circ(f_0,\ldots,f_{q-1},r,r+1)$. In all
cases the coefficient of~$i$ in $u\circ f$ is nonnegative.

The remaining case, $v\circ f$ when the image of~$f$ contains~$r$
but not $r+1$, is equivalent to the previous case, because $v\circ
f=u\circ f'$, where $f'$~is got from~$f$ by changing~$r$ to $r+1$.

This completes the proof.
\end{proof}

\begin{proposition} \label{factor middle}
Let $x$ be a morphism in $\Or(m,n)$ with $m>0$ and $x[m]=[t]$ such
that $x$~is a linear combination of terms $(a_0,\ldots,a_m)$ with
$a_{m-1}<t$. Then there is a factorisation
$$x=u\vee_{m-1}v$$
in~$\Or$ such that $u[m]=[t']$ with $t'<t$ and $v[m]=[t]$ and
$v$~is a linear combination of terms $(a_0,\ldots,a_m)$ with
$a_{m-1}=a_m$ or $a_m=t$.
\end{proposition}

\begin{proof} We get a factorisation $x=u\vee_{m-1}v$ in~$\Or$
just as in the proof of Proposition~\ref{factor filler start},
except that $r$~is replaced by $m-1$, and we find that $u$~and~$v$
are as described.
\end{proof}

\begin{proposition} \label{factor filler finish}
Let $x$ be a morphism in $\Or(m,n)$ with $x[m]=[t]$, and let $r$
be an integer with $0\leq r\leq m-2$ such that $x$~is a linear
combination of terms $(a_0,\ldots,a_m)$ with $a_{r+1}=a_m$ or
$a_m=t$. Then there is a factorisation
$$x=u\vee_r v=(\partial_{r+2}u\tri_r\partial_r u)\vee_r v$$
in~$\Or$ such that $v[m]=[t]$ and $v$~is a linear combination of
terms $(a_0,\ldots,a_m)$ with $a_r=a_m$ or $a_m=t$.
\end{proposition}

\begin{proof} Let $\alpha,\beta\colon\Z\Delta(m,n)\to\Z\Delta(m,n)$ be the
linear functions given on morphisms $(a_0,\ldots,a_m)$ in
$\Delta(m,n)$ as follows:
\begin{align*}
&\alpha(a_0,\ldots,a_m)=
 \left\{\begin{alignedat}{2}
 &(a_0,\ldots,a_{r-1},a_r,a_{r+1},a_{r+2}\ldots,a_m)
 &\quad &\text{if $a_m<t$,}\\
 &(a_0,\ldots,a_{r-1},a_r,a_r,a_{r+2}\ldots,a_m)
 &&\text{if $a_m\geq t$,}
 \end{alignedat}\right.\\
&\beta(a_0,\ldots,a_m)=
 \left\{\begin{alignedat}{2}
 &(a_0,\ldots,a_{r-1},a_{r+1},a_{r+1},a_{r+2}\ldots,a_m)
 &\quad &\text{if $a_m<t$,}\\
 &(a_0,\ldots,a_{r-1},a_r,a_{r+1},a_{r+2}\ldots,a_m)
 &&\text{if $a_m\geq t$}
 \end{alignedat}\right.
\end{align*}
(this is the same as in the proof of Proposition~\ref{factor
filler start}, except that the splitting depends on~$a_m$ rather
than~$a_{r+1}$). We now get a factorisation
$$x=u\vee_r v=(\partial_{r+2}u\tri_r\partial_r u)\vee_r v$$
in $\Z\Delta(m,n)$, where $u=\alpha x$ and $v=\beta x$. It is
clear that $v[m]=x[m]$, so that $v[m]=[t]$, and it clear that
$v$~is a linear combination of terms $(a_0,\ldots,a_m)$ with
$a_r=a_m$ or $a_m=t$. It therefore remains to show that
$u$~and~$v$ are in $\Or(m,n)$. Obviously the sums of the
coefficients in $u$~and~$v$ are equal to~$1$, so it suffices to
show that an injective morphism $i=(i_0,\ldots,i_q)$ with
codomain~$n$ has a nonnegative coefficient in $u\circ f$ and
$v\circ f$ whenever $f=(f_0,\ldots,f_q)$ is an injective morphism
with codomain~$m$. The case $u\circ f$ with $r+1$ not in the image
of~$f$ and the case $v\circ f$ with $r$ not in the image of~$f$
are handled as in the proof of Proposition~\ref{factor filler
start}. Also as in the proof of that proposition, the case $u\circ
f$ with the image of~$f$ containing $r+1$ but not~$r$ reduces to
the case $v\circ f$ with the image of~$f$ containing~$r$ but
not~$r+1$. In the remaining cases, we argue as follows.

Suppose that the image of~$f$ contains $r$~and~$r+1$. If
$f_q>r+1$, then the coefficient of~$i$ in $u\circ f$ is zero and
the coefficient of~$i$ in $v\circ f$ is as in $x\circ f$. If
$f_{q-1}=r$ and $f_q=r+1$ and $i_q<t$ then the coefficient of~$i$
in $u\circ f$ is the coefficient of~$i$ in
$x\circ(f_0,\ldots,f_{q-1},m)$ and the coefficient of~$i$ in
$v\circ f$ is the coefficient of $(i_0,\ldots,i_q,t)$ in
$x\circ(f_0,\ldots,f_q,m)$. If $f_{q-1}=r$ and $f_q=r+1$  and
$i_q\geq t$ then the coefficient of~$i$ in $u\circ f$ is zero and
the coefficient of~$i$ in $v\circ f$ is as in $x\circ f$.

Now suppose that the image of~$f$ contains~$r$ but not $r+1$, say
$f_j=r$. If $i_q\geq t$, then the coefficient of~$i$ in $v\circ f$
is the same as in $x\circ f$. If $i_q<t$ and $j<q$ and $f_q=m$,
then the coefficient of~$i$ in $v\circ f$ is zero. If $i_q<t$ and
$j<q$ and $f_q<m$, then the coefficient of~$i$ in $v\circ f$ is
the coefficient of $(i_0,\ldots,i_q,t)$ in
$x\circ(f_0,\ldots,f_q,m)$. Finally, if $i_q<t$ and $j=q$, then
the coefficient of~$i$ in $v\circ f$ is the coefficient of~$i$ in
$x\circ(f_0,\ldots,f_{q-1},m)$ plus the coefficient of
$(i_0,\ldots,i_q,t)$ in $x\circ(f_0,\ldots,f_q,m)$. This completes
the proof.
\end{proof}

\begin{proposition} \label{constant end}
Let $x$ be a morphism in $\Or(m,n)$ with $x[m]=[t]$ such that
$x$~is a linear combination of terms $(a_0,\ldots,a_m)$ with
$a_0=a_m$ or $a_m=t$. Then $x$~is a linear combination of terms
$(a_0,\ldots,a_m)$ with $a_m=t$.
\end{proposition}

\begin{proof} For $a<t$, we must show that $(a,\ldots,a)$ has zero
coefficient in~$x$. But this holds because the coefficient of
$(a,\ldots,a)$ in~$x$ is the coefficient of~$[a]$ in $x[m]$, and
this coefficient is zero by hypothesis.
\end{proof}

We now combine these results to get factorisations for arbitrary
morphisms in~$\Or$.

\begin{theorem} \label{factorisations} Let $x$ be a morphism
in~$\Or$. Then $x$~can be factorised into morphisms in~$\Delta$ by
using the operations $\tri_i$~and~$\vee_i$.
\end{theorem}

\begin{proof} By Proposition~\ref{first and last}, there is an
integer~$t$ with $0\leq t\leq n$ such that $x[m]=[t]$, and $t$~is
the largest integer to appear in the terms of~$x$. We will use
induction on~$m$, and for a fixed~$m$ we will use induction
on~$t$.

Suppose first that $x$~has a term $(t,\ldots,t)$ (this includes
the base cases with $m=0$ or $t=0$). By
Proposition~\ref{constant}, $x$~is equal to $(t,\ldots,t)$, so
$x$~is already a morphism in~$\Delta$.

From now on, suppose that $x$~has no term $(t,\ldots,t)$. It
follows that $x$~is a linear combination of terms
$(a_0,\ldots,a_m)$ with $a_0<t$. We now apply Propositions
\ref{factor filler start}--\ref{factor filler finish}. By applying
Proposition~\ref{factor filler start} for $r=0$, $r=1$, \dots,
$r=m-2$ successively, we get a factorisation
$$x=
 [\,\ldots[\tilde x
 \vee_{m-2}(v'_{m-2}\tri_{m-2}v''_{m-2})]\vee_{m-3}
 \ldots\,]
 \vee_0(v'_0\tri_0 v''_0)$$
in~$\Or$ such that $\tilde x[m]=[t]$ and $\tilde x$~is a linear
combination of terms $(a_0,\ldots,a_m)$ with $a_{m-1}<t$. By
applying Proposition~\ref{factor middle} we get a factorisation
$$\tilde x=y\vee_{m-1}\tilde y$$
in~$\Or$ such that $y[m]=[t']$ with $t'<t$ and $\tilde y[m]=[t]$
and $\tilde y$~is a linear combination of terms $(a_0,\ldots,a_m)$
with $a_{m-1}=a_m$ or $a_m=t$. Finally, by applying
Proposition~\ref{factor filler finish} for $r=m-2$, $r=m-3$,
\dots, $r=0$ successively, we get a factorisation
$$\tilde y=
 (u'_{m-2}\tri_{m-2}u''_{m-2})\vee_{m-2}[\,\ldots
 \vee_1[(u'_0\tri_0 u''_0)\vee_0 z]\ldots\,]$$
in~$\Or$ such that $z[m]=[t]$ and $z$~is a linear combination of
terms $(a_0,\ldots,a_m)$ with $a_0=a_m$ or $a_m=t$.

We have now got a factorisation of~$x$ in~$\Or$ with factors
$v'_r$, $v''_r$, $y$, $u'_r$, $u''_r$ and~$z$. The factors $v'_r$,
$v''_r$, $u'_r$ and $u''_r$ are in $\Or(m-1,n)$, so they can be
factorised into morphisms in~$\Delta$ by the induction on~$m$. The
factor~$y$ is in $\Or(m,n)$ with $y[m]=[t']$ for some $t'<t$, so
it can be factorised into morphisms in~$\Delta$ by the induction
on~$t$. By Proposition~\ref{constant end}, the factor~$z$ is a
linear combination of terms $(a_0,\ldots,a_m)$ with $a_m=t$, so it
can be expressed as $(z_t,t)$ in the notation of
Proposition~\ref{tails}. Here $z_t=\partial_m z$, so $z_t$~is in
$\Or(m-1,n)$. By the induction on~$m$, we can factorise~$z_t$ into
morphisms in~$\Delta$, and this clearly induces a similar
factorisation for~$z$. We have now factorised all the factors
of~$x$, so we have a factorisation of~$x$ itself. This completes
the proof.
\end{proof}

For example, suppose that $x\in\Or(1,n)$. Then the factorisation
reduces to repeated application of Proposition~\ref{factor
middle}, and it takes the form
$$x=(i_0,i_0)\vee_0(i_0,i_1)\vee_0\ldots\vee_0(i_{q-1},i_q),$$
where $q\geq 0$ and $0\leq i_0<i_1<\ldots<i_q\leq n$; in other
words, $x=(i_0,i_0)$ or
$x=(i_0,i_1)\vee_0\ldots\vee_0(i_{q-1},i_q)$ as in
Example~\ref{filler examples}.

Finally, from Theorem~\ref{factorisations}, we get the main
result.

\begin{theorem} \label{decompositions} The category of
orientals~$\Or$ is the closure of the subcategory of
simplexes~$\Delta$ under the filler operations and composition.
\end{theorem}

\begin{proof} By Proposition~\ref{closure}, $\Or$~is closed under
the filler and pasting operations. Theorem~\ref{factorisations}
therefore shows that $\Or$~is the closure of~$\Delta$ under the
filler and pasting operations. If $x\vee_i y$ is a pasting, then
$x\vee_i y=\partial_{i+1}(x\tri_i y)$ by Proposition~\ref{faces},
so $x\vee_i y$ has the form $(x\tri_i y)\circ\delta$ with $\delta$
a morphism in~$\Delta$. Therefore $\Or$~is the closure of~$\Delta$
under fillers and composition. This completes the proof.
\end{proof}


\providecommand{\bysame}{\leavevmode\hbox
to3em{\hrulefill}\thinspace}
\providecommand{\MR}{\relax\ifhmode\unskip\space\fi MR }
\providecommand{\MRhref}[2]{%
  \href{http://www.ams.org/mathscinet-getitem?mr=#1}{#2}
} \providecommand{\href}[2]{#2}

\end{document}